\documentclass[12pt]{article}
\usepackage{amssymb}

\newcommand{\be}{\begin{equation}}
\newcommand{\ee}{\end{equation}}
\newcommand{\bea}{\begin{eqnarray}}
\newcommand{\eea}{\end{eqnarray}}
\newcommand{\barray}{\begin{array}}
\newcommand{\earray}{\end{array}}
\newcommand{\pa}{\partial}
\newcommand{\nn}{\nonumber}
\newcommand{\bitem}{\begin{itemize}}
\newcommand{\eitem}{\end{itemize}}
\newtheorem{teo}{Theorem}[section]
\newcommand{\bt}{\begin{teo}}
\newcommand{\et}{\end{teo}}
\newtheorem{Def}{Definition}[section]
\newcommand{\bd}{\begin{Def}}
\newcommand{\ed}{\end{Def}}
\newtheorem{lem}{Lemma}[section]
\newcommand{\bl}{\begin{lem}}
\newcommand{\el}{\end{lem}}
\newtheorem{prop}{Proposition}[section]
\newcommand{\bp}{\begin{prop}}
\newcommand{\ep}{\end{prop}}
\newtheorem{cor}{Corollary}[section]
\newcommand{\bc}{\begin{cor}}
\newcommand{\ec}{\end{cor}}
\newtheorem{ex}{Example}[section]
\newcommand{\bex}{\begin{ex}}
\newcommand{\eex}{\end{ex}}
\newtheorem{rem}{Remark}[section]
\newcommand{\br}{\begin{rem}}
\newcommand{\er}{\end{rem}}

\begin{document}

\begin{center}
{\Large \textbf{On compatible metrics and diagonalizability \\ of
non-locally bi-Hamiltonian systems \\ of hydrodynamic
type\footnote{The work was supported by the Max-Planck-Institut
f\"{u}r Mathematik (Bonn, Germany), by the Russian Foundation for
Basic Research (project no.~08-01-00464) and by the programme
``Leading Scientific Schools'' (project no. NSh-1824.2008.1).}}}
\end{center}

\medskip

\begin{center}
{\large \bf {O. I. Mokhov}}
\end{center}

\bigskip

\begin{center}
\bf {Abstract}
\end{center}

We study bi-Hamiltonian systems of hydrodynamic type with
non-singular (semi\-simple) non-local bi-Hamiltonian structures and
prove that such systems of hydrodynamic type are diagonalizable.
Moreover, we prove that for an arbitrary non-singular (semisimple)
non-locally bi-Hamiltonian system of hydrodynamic type, there exist
local coordinates (Riemann invariants) such that all the related
matrix differential-geometric objects, namely, the matrix $V^i_j
(u)$ of this system of hydrodynamic type, the metrics $g^{ij}_1 (u)$
and $g^{ij}_2 (u)$ and the affinors $(w_{1, n})^i_j (u)$ and $(w_{2,
n})^i_j (u)$ of the non-singular non-local bi-Hamiltonian structure
of this system, are diagonal in these local coordinates. The proof
is a natural consequence of the general results of the theory of
compatible metrics and the theory of non-local bi-Hamiltonian
structures developed earlier by the present author in [21]--[33].

\bigskip

\begin{center}
{\large \bf {Introduction}}
\end{center}

\medskip

In this paper we consider $(1+1)$-dimensional non-singular
(semisimple) non-locally bi-Hamiltonian systems of hydrodynamic type
and prove their diagonalizability. Moreover, we prove that for an
arbitrary non-singular (semisimple) non-locally bi-Hamiltonian
system of hydrodynamic type, there exist local coordinates (Riemann
invariants) such that all the related matrix differential-geometric
objects, namely, the matrix $V^i_j (u)$ of this system of
hydrodynamic type, the metrics $g^{ij}_1 (u)$ and $g^{ij}_2 (u)$ and
the affinors $(w_{1, n})^i_j (u)$ and $(w_{2, n})^i_j (u)$ of the
non-singular non-local bi-Hamiltonian structure of this system, are
diagonal in these local coordinates. Let us give here very briefly
basic well-known notions and results necessary for us. Recall that
$(1+1)$-dimensional {\it systems of hydrodynamic type\,} [1] are
arbitrary $(1+1)$-dimensional evolution quasilinear systems of
first-order partial differential equations, i.e., equations of the
form \be u^i_t = V^i_j (u) u^j_x, \ \ \ \ 1 \leq i, j \leq N,
\label{1} \ee where $u = (u^1, \ldots, u^N)$ are local coordinates
on a certain smooth $N$-dimensional manifold or in a domain of
${\mathbb{R}}^N$ (or ${\mathbb{C}}^N$); $u^i (x)$ are functions
(fields) of one variable $x$ that are evolving with respect to $t$;
$V^i_j (u)$ is an arbitrary $(N \times N)$-matrix depending on $u$
(this matrix is a mixed tensor of the type (1,\,1), i.e., an
affinor, with respect to local changes of coordinates $u$).

\begin{center}
{\large \bf {Non-locally Hamiltonian systems of hydrodynamic type}}
\end{center}

We will consider systems of the form (\ref{1}) that are Hamiltonian
with respect to arbitrary non-degenerate non-local Poisson brackets
of hydrodynamic type (the Ferapontov brackets [2], see also [3] and
[1] for the Mokhov--Ferapontov brackets and the Dubrovin--Novikov
brackets in partial non-local and local cases, respectively), i.e.,
\be u^i_t = V^i_j (u) u^j_x = \{ u^i (x), H\}, \ \ \ \ 1 \leq i, j
\leq N, \label{nonlmS1}\ee where the functional \be H = \int h
(u(x)) dx \label{nonlmS2a}\ee is the Hamiltonian of the system
(\ref{nonlmS1}) (the function $h (u)$ is the density of the
Hamiltonian) and the Poisson bracket has the form \be \{ u^i (x),
u^j (y)\} = P^{ij} \delta (x - y), \ \ \ \ 1 \leq i, j \leq N,
\label{nonlm1} \ee \bea && P^{ij} =
g^{ij}(u(x)) {d \over dx} + b^{ij}_k (u(x))\, u^k_x + \nn\\
&& + \sum_{m, n =1}^L \mu^{mn} (w_m)^i_k (u (x)) u^k_x \left ( {d
\over dx} \right )^{-1} \circ (w_n)^j_s (u (x)) u^s_x,
\label{nonlm2} \eea where the coefficients $g^{ij} (u),$ $b^{ij}_k
(u),$ and $(w_n)^i_j (u),$ $1 \leq i, j, k \leq N,$ $1 \leq n \leq
L,$ are smooth functions of local coordinates, $\det (g^{ij} (u))
\neq 0,$ $\mu^{mn}$ is an arbitrary non-degenerate symmetric
constant matrix, $\mu^{mn} = \mu^{nm},$ $\mu^{mn} = {\rm const},$
$\det (\mu^{mn})\neq0$. For two arbitrary functionals $I$ and $J$
the Poisson bracket (\ref{nonlm1}), (\ref{nonlm2}) has the form \be
\{ I, J\} = \int {\delta I \over \delta u^i (x)} P^{ij} {\delta J
\over \delta u^j (x)} dx. \label{nonlm3} \ee Poisson brackets of the
form (\ref{nonlm2}), (\ref{nonlm3}) were introduced and studied by
Ferapontov in [2]; these brackets are a non-local generalization of
the Dubrovin--Novikov brackets (local Poisson brackets of
hydrodynamic type generated by flat metrics $g^{ij} (u)$; there are
no non-local terms in this case, $L = 0$, or $(w_n)^i_j (u) = 0$)
[1] and the Mokhov--Ferapontov brackets (non-local Poisson brackets
of hydrodynamic type generated by metrics of constant curvature $K$;
in this case $L = 1$, $\mu^{11} = K$, $(w_1)^i_j (u) = \delta^i_j$)
[3]. Ferapontov proved that a non-local operator $P^{ij}$ of the
form (\ref{nonlm2}) gives a Poisson bracket (\ref{nonlm3}) if and
only if there is an $N$-dimensional submanifold with flat normal
bundle in an $(N+L)$-dimensional pseudo-Euclidean space such that
$g^{ij} (u)$ is the contravariant first fundamental form; $b^{ij}_k
(u) = - g^{is} (u) \Gamma^j_{sk} (u)$; $\Gamma^j_{sk} (u)$ are the
Christoffel symbols of the Levi-Civita connection of the metric
$g^{ij} (u)$; $(w_n)^i_j (u),$ $1 \leq n \leq L,$ are the Weingarten
operators (the Weingarten affinors) of the submanifold; and
$\mu^{mn}$ is the Gram matrix of the corresponding parallel bases in
the normal spaces of the submanifold (all torsion forms of the
submanifold with flat normal bundle vanish in these bases in the
normal spaces).

In other words, the non-local operator (\ref{nonlm2}) gives a
Poisson bracket (\ref{nonlm3}) if and only if its coefficients
satisfy the relations (see also [4]) \be g^{ij} = g^{ji}, \label{01}
\ee \be {\pa g^{ij} \over \pa u^k} = b^{ij}_k + b^{ji}_k, \label{02}
\ee \be g^{is} b^{jk}_s = g^{js} b^{ik}_s, \label{03} \ee \be g^{is}
(w_n)^j_s = g^{js} (w_n)^i_s, \label{04} \ee \be (w_n)^i_s (w_m)^s_j
= (w_m)^i_s (w_n)^s_j, \label{05} \ee \be g^{is} g^{jr} {\pa
(w_n)^k_r \over \pa u^s} - g^{jr} b^{ik}_s (w_n)^s_r  = g^{js}
g^{ir} {\pa (w_n)^k_r \over \pa u^s} - g^{ir} b^{jk}_s (w_n)^s_r,
\label{06} \ee \be g^{is} \left ( {\pa b^{jk}_s \over \pa u^r} -
{\pa b^{jk}_r \over \pa u^s} \right ) + b^{ij}_s b^{sk}_r - b^{ik}_s
b^{sj}_r = \sum_{m = 1}^L \sum_{n =1}^L \mu^{mn} g^{is} \left (
(w_m)^j_r (w_n)^k_s -
 (w_m)^j_s (w_n)^k_r \right ). \label{07}
\ee

\begin{center}
{\large \bf {Non-locally Hamiltonian affinors}}
\end{center}

The Hamiltonian $H$ of the system (\ref{nonlmS1})--(\ref{nonlm3})
must also be a first integral of all the systems of hydrodynamic
type that are given by the affinors $(w_n)^i_j (u),$ $1 \leq n \leq
L,$ of the non-local operator (\ref{nonlm2}) (these systems are
called the {\it structural flows} of the non-local Poisson bracket
(\ref{nonlm1})--(\ref{nonlm3})) [2]: \be u^i_{t_n} = (w_n)^i_j (u)
u^j_x, \ \ \ H_{t_n} = 0, \ \ \ 1 \leq n \leq L. \label{strhy} \ee
For each $n$, $1 \leq n \leq L$, there exist a function $f_n (u)$
such that \be {\pa h \over \pa u^j} (w_n)^j_s (u) = {\pa f_n \over
\pa u^s}, \ \ \ 1 \leq n \leq L. \label{strhy2}\ee In this case the
affinor $V^i_j (u)$ of the system of hydrodynamic type
(\ref{nonlmS1})--(\ref{nonlm3}) has the form \be V^i_j (u) =
g^{is}(u) {\pa^2 h \over \pa u^s \pa u^j} - g^{is} (u) \Gamma^p_{sj}
(u) {\pa h \over \pa u^p} + \sum_{m,  n =1}^L \mu^{mn} (w_m)^i_j (u)
f_n (u), \label{nonlmp2} \ee i.e., \bea  && V^i_j (u) = g^{is}(u)
\nabla_s \nabla_j h (u)  + \sum_{m, n =1}^L \mu^{mn} (w_m)^i_j (u)
f_n (u) = \nn \\ &&= \nabla^i \nabla_j h (u)  + \sum_{m, n =1}^L
\mu^{mn} (w_m)^i_j (u) f_n (u), \label{nonlmp3} \eea where
$\nabla_k$ is the covariant differentiation generated by the
Levi-Civita connection $\Gamma^j_{sk} (u)$ of the metric $g^{ij}
(u).$

We will call an affinor $V^i_j (u)$ {\it Hamiltonian} (or {\it
non-locally Hamiltonian}) if there exist an $N$-dimensional
submanifold with flat normal bundle in an $(N+L)$-dimensional
pseudo-Euclidean space and functions $h (u)$ and $f_n (u)$, $1 \leq
n \leq L$, such that the affinor $V^i_j (u)$ has the form
(\ref{nonlmp3}), where $g^{ij} (u)$ is the contravariant first
fundamental form of the submanifold; $\Gamma^j_{sk} (u)$ are the
Christoffel symbols of the Levi-Civita connection of the metric
$g^{ij} (u)$; $(w_n)^i_j (u),$ $1 \leq n \leq L,$ are the Weingarten
operators of the submanifold; and $\mu^{mn}$ is the Gram matrix of
the corresponding parallel bases in the normal spaces of the
submanifold (such that all torsion forms of the submanifold with
flat normal bundle vanish in these bases in the normal spaces), and
the functions $h (u)$ and $f_n (u)$, $1 \leq n \leq L$, satisfy
relations (\ref{strhy2}): \be (w_n)^j_s (u) \nabla_j h (u) =
\nabla_s f_n (u), \ \ \ 1 \leq n \leq L. \label{strhy3}\ee In this
case we will also speak that the affinor $V^i_j (u)$
(\ref{nonlmp3}), (\ref{strhy3}) is Hamiltonian with respect to the
corresponding non-local Poisson bracket of hydrodynamic type
(\ref{nonlm1})--(\ref{nonlm3}). Obviously, this definition is
invariant. We note that one can also consider it as a definition of
a non-locally Hamiltonian system of hydrodynamic type (\ref{1}).

Affinors that are Hamiltonian with respect to the Dubrovin--Novikov
brackets ({\it locally Hamiltonian affinors}) were studied in detail
by Tsarev in the remarkable work [5]. Affinors that are Hamiltonian
with respect to the non-local Mokhov--Ferapontov brackets (the
Mokhov--Ferapontov affinors) were studied in detail in the paper
[3].

Using relations (\ref{01})--(\ref{07}) it is easy to prove that the
following relations always hold for non-locally Hamiltonian affinors
$V^i_j (u)$: \be g_{is} (u) V^s_j (u) = g_{js} (u) V^s_i
(u),\label{jj1}\ee \be \nabla_j V^i_k (u) = \nabla_k V^i_j (u).
\label{jj2}\ee Here $g_{ij} (u)$ is the inverse of the matrix
$g^{ij} (u)$, $g_{is} (u) g^{sj} (u) = \delta^j_i$ (the covariant
metric).

For locally Hamiltonian affinors these important relations are very
simple in flat local coordinates of the metric $g_{ij} (u)$ and
Tsarev proved that in this flat case relations (\ref{jj1}) and
(\ref{jj2}) are not only necessary but also sufficient for an
affinor to be locally Hamiltonian (an affinor $V^i_j (u)$ is locally
Hamiltonian if and only if there exists a flat metric $g_{ij} (u)$
such that relations (\ref{jj1}) and (\ref{jj2}) hold) [5]. This
result was generalized to the case of the Mokhov--Ferapontov
affinors in [3]: an affinor $V^i_j (u)$ is Hamiltonian with respect
to a non-local Mokhov--Ferapontov bracket if and only if there
exists a metric $g_{ij} (u)$ of constant curvature such that
relations (\ref{jj1}) and (\ref{jj2}) hold.

Let the affinor $V^i_j (u)$ of a system of hydrodynamic type
(\ref{1}) satisfy relations (\ref{jj1}) and (\ref{jj2}), i.e., there
exists a metric $g_{ij} (u)$ such that relations (\ref{jj1}) and
(\ref{jj2}) hold. If this system of hydrodynamic type is
diagonalizable, i.e., there exist local coordinates such that the
affinor $V^i_j (u)$ is a diagonal matrix $V^i_j (u) = V^i (u)
\delta^i_j$ in these special local coordinates (such local
coordinates are called {\it Riemann invariants}), and strictly
hyperbolic, i.e., all the eigenvalues $V^i (u)$, $1 \leq i \leq N$,
are distinct ($V^i (u) \neq V^j (u)$ when $i \neq j$), then it can
be integrated by the generalized hodograph method (Tsarev, see [5]).
In this case, relation (\ref{jj1}) is equivalent to the condition
that the metric $g_{ij} (u)$ is also diagonal in these special local
coordinates, $g_{ij} (u) = g_i (u) \delta_{ij}$, i.e., the Riemann
invariants are orthogonal curvilinear coordinates in the
corresponding pseudo-Riemannian space, and relation (\ref{jj2}) is
equivalent to the condition \be {\pa V^i \over \pa u^j} = {\pa \ln
{\sqrt{g_i (u)}} \over \pa u^j} (V^j (u) - V^i (u)), \ \ \ i \neq
j.\ee Hence, the following relation holds for the eigenvalues $V^i
(u)$: \be {\pa \over \pa u^k} \left ( {1 \over (V^j (u) - V^i (u))}
{\pa V^i \over \pa u^j} \right ) = {\pa \over \pa u^j} \left ( {1
\over (V^k (u) - V^i (u))} {\pa V^i \over \pa u^k} \right ).
\label{semiham} \ee A strictly hyperbolic diagonal system of
hydrodynamic type is called {\it semi-Hamiltonian} if relations
(\ref{semiham}) hold (Tsarev, see [5]). In [5] Tsarev proved that
any strictly hyperbolic diagonal semi-Hamiltonian system of
hydrodynamic type is integrable by the generalized hodograph method.

\begin{center}
{\large \bf {Diagonalizable affinors}}
\end{center}

Recall that the very important problem of diagonalizability for an
affinor, which had been posed, in fact, by Riemann, was completely
solved by Haantjes in [6] on the base of the previous Nijenhuis'
results [7]. An affinor $V^i_j (u)$ is diagonalizable by a local
change of coordinates in a domain if and only if it is
diagonalizable at any point and its Haantjes tensor vanishes. The
Haantjes tensor of an affinor $V (u) = V^i_j (u)$ is the following
tensor of the type $(1, 2)$ (a skew-symmetric vector-valued 2-form)
generated by the affinor $V^i_j (u)$: \be H(X, Y) = N(V(X), V(Y)) +
V^2(N(X, Y)) - V(N(X, V(Y))) - V(N(V(X), Y)),\ee where $X (u)$ and
$Y (u)$ are arbitrary vector fields, $V(X)$ is the vector field
$V^i_j (u) X^j (u)$, $N(X, Y)$ is the Nijenhuis tensor of the
affinor $V^i_j (u)$, i.e., the following tensor of the type $(1, 2)$
(a skew-symmetric vector-valued 2-form) generated by the affinor
$V^i_j (u)$: \be N(X, Y) = [V(X), V(Y)] + V^2([X, Y]) - V([X, V(Y)])
- V([V(X), Y]),\ee where $[X, Y]$ is the commutator of the vector
fields $X (u)$ and $Y (u)$. In components, the Nijenhuis tensor of
the affinor $V^i_j (u)$ has the form \be N^k_{ij} (u) = V^s_i (u)
{\pa V^k_j \over \pa u^s} - V^s_j (u) {\pa V^k_i \over \pa u^s}  +
V^k_s (u) {\pa V^s_i \over \pa u^j} - V^k_s (u) {\pa V^s_j \over \pa
u^i}  \label{nij01} \ee and the Haantjes tensor of the affinor
$V^i_j (u)$ has the form \bea && H^i_{jk} (u) = V^i_s (u) V^s_r (u)
N^r_{jk} (u) - V^i_s (u) N^s_{rk} (u) V^r_j (u) - \nn\\ && - V^i_s
(u) N^s_{jr} (u) V^r_k (u) + N^i_{sr} (u) V^s_j (u) V^r_k (u). \eea
Recall also that invariant tensor conditions that a strictly
hyperbolic system of hydrodynamic type is semi-Hamiltonian were
found in [8].

\begin{center}
{\large \bf {Non-locally bi-Hamiltonian systems of hydrodynamic
type}}
\end{center}

We will consider bi-Hamiltonian systems of hydrodynamic type. Recall
that two Poisson brackets are called {\it compatible} if any linear
combination of these Poisson brackets is also a Poisson bracket [9],
and a system of equations that is Hamiltonian with respect to two
linearly independent compatible Poisson brackets is called {\it
bi-Hamiltonian}. In this paper we will consider systems of
hydrodynamic type that are bi-Hamiltonian with respect to two
linearly independent compatible non-degenerate non-local Poisson
brackets of hydrodynamic type (\ref{nonlm1})--(\ref{nonlm3}), \be
u^i_t = V^i_j (u) u^j_x = \{ u^i (x), H_1\}_1 = \{ u^i (x), H_2\}_2,
\ \ \ \ 1 \leq i, j \leq N, \label{nonlmS2}\ee \be H_1 = \int h_1
(u(x)) dx, \ \ \ H_2 = \int h_2 (u(x)) dx, \label{nonlm3b} \ee \be
\{ I, J\}_1 = \int {\delta I \over \delta u^i (x)} P_1^{ij} {\delta
J \over \delta u^j (x)} dx, \ \ \ \{ I, J\}_2 = \int {\delta I \over
\delta u^i (x)} P_2^{ij} {\delta J \over \delta u^j (x)} dx,
\label{nonlm3a}\ee \bea && P_1^{ij} = g_1^{ij}(u(x)) {d \over dx} +
b^{ij}_{1, k} (u(x))\, u^k_x + \nn \\ && + \sum_{m, n =1}^L
\mu_1^{mn} (w_{1, m})^i_k (u (x)) u^k_x \left ( {d \over dx} \right
)^{-1} \circ (w_{1, n})^j_s (u (x)) u^s_x, \label{nonlm4} \eea \bea
&& P_2^{ij} = g_2^{ij}(u(x)) {d \over dx} + b^{ij}_{2, k} (u(x))\,
u^k_x + \nn \\ && + \sum_{m, n =1}^L \mu_2^{mn} (w_{2, m})^i_k (u
(x)) u^k_x \left ( {d \over dx} \right )^{-1} \circ (w_{2, n})^j_s
(u (x)) u^s_x. \label{nonlm5} \eea

We will call an affinor $V^i_j (u)$ {\it bi-Hamiltonian} (or {\it
non-locally bi-Hamiltonian}) if this affinor is Hamiltonian with
respect to two linearly independent compatible non-degene\-rate
non-local Poisson brackets of hydrodynamic type
(\ref{nonlm1})--(\ref{nonlm3}). This definition is invariant.

In this paper we prove that $(1+1)$-dimensional non-singular
(semisimple) non-locally bi-Hamiltonian systems of hydrodynamic type
are diagonalizable. Recall that a pair of pseudo-Riemannian metrics
$g_1^{ij} (u)$ and $g_2^{ij} (u)$ is called {\it non-singular} (or
{\it semisimple}) if the eigenvalues of this pair of metrics, i.e.,
the roots of the equation \be \det ( g_1^{ij} (u) -  \lambda
g_2^{ij} (u)) =0,\ee are distinct. In this case the non-locally
bi-Hamiltonian system of hydrodynamic type
(\ref{nonlmS2})--(\ref{nonlm5}) and the corresponding bi-Hamiltonian
affinor are also called {\it non-singular} (or {\it semisimple}).

It is important to note that, generally speaking, integrable
bi-Hamiltonian systems of hydrodynamic type are not necessarily
diagonalizable if we consider an other class of compatible Poisson
brackets (even if both the compatible Poisson brackets are local).
This is a nontrivial fact and we give here a very important example
in detail.

{\bf Example} ({\it a non-diagonalizable integrable bi-Hamiltonian
system of hydrodynamic type {\rm [10]-[12]}}). Let us consider the
associativity equations of two-dimensional topological quantum field
theories (the Witten--Dijkgraaf--Verlinde--Verlinde equations, see
[13]--[16]) for a function (a {\it potential\/}) $\Phi = \Phi (u^1,
\ldots, u^N)$, \be \sum_{k = 1}^N \sum_{l = 1}^N {\pa^3 \Phi \over
\pa u^i \pa u^j \pa u^k} \eta^{kl} {\pa^3 \Phi \over \pa u^l \pa u^m
\pa u^n} = \sum_{k = 1}^N \sum_{l = 1}^N {\pa^3 \Phi \over \pa u^i
\pa u^m \pa u^k} \eta^{kl} {\pa^3 \Phi \over \pa u^l \pa u^j \pa
u^n}, \label{ass1} \ee where $\eta^{ij}$ is an arbitrary constant
nondegenerate symmetric matrix, $\eta^{ij} = \eta^{ji},$ $\eta^{ij}
= {\rm const},$ $\det (\eta^{ij})\neq0$. We recall that the
associativity equations (\ref{ass1}) are consistent and integrable
by the inverse scattering method, they possess a rich set of
nontrivial solutions, and each solution $\Phi (u^1, \ldots, u^N)$ of
the associativity equations (\ref{ass1}) gives $N$-parameter
deformations of special Frobenius algebras (some special commutative
associative algebras equipped with nondegenerate invariant symmetric
bilinear forms) (see [13]). Indeed, consider algebras $A (u)$ in an
$N$-dimensional vector space with the basis $e_1, \ldots, e_N$ and
the multiplication (see [13]) \be e_i \circ e_j = c^k_{ij} (u) e_k,
\ \ \ \ c^k_{ij} (u) = \eta^{ks} {\pa^3 \Phi \over \pa u^s \pa u^i
\pa u^j}. \label{al1} \ee For all values of the parameters $u =
(u^1, \ldots, u^N)$ the algebras $A (u)$ are commutative, $e_i \circ
e_j = e_j \circ e_i,$ and the associativity condition \be (e_i \circ
e_j) \circ e_k = e_i \circ (e_j \circ e_k) \label{al2} \ee in the
algebras $A (u)$ is equivalent to equations (\ref{ass1}). The matrix
$\eta_{ij}$ inverse to the matrix $\eta^{ij}$, $\eta^{is} \eta_{sj}
= \delta^i_j$, defines a nondegenerate invariant symmetric bilinear
form on the algebras $A (u)$, \be \langle e_i, e_j \rangle =
\eta_{ij}, \ \ \ \ \langle e_i \circ e_j, e_k \rangle = \langle e_i,
e_j \circ e_k \rangle. \label{al3} \ee Recall that locally the
tangent space at every point of any Frobenius manifold (see [13])
possesses the structure of Frobenius algebra
(\ref{al1})--(\ref{al3}), which is determined by a solution of the
associativity equations (\ref{ass1}) and smoothly depends on the
point. Let $N = 3$ and the metric $\eta_{ij}$ be antidiagonal \be
(\eta_{ij}) =
\left ( \begin{array} {ccc} 0&0&1\\
0&1&0\\
1&0&0
\end{array} \right ),
\ee and the function $\Phi (u)$ has the form
$$
\Phi (u) = {1 \over 2} (u^1)^2 u^3 + {1 \over 2} u^1 (u^2)^2 + f
(u^2, u^3).
$$
In this case $e_1$ is the unit in the Frobenius algebra
(\ref{al1})--(\ref{al3}), and the associativity equations
(\ref{ass1}) for the function $\Phi (u)$ are equivalent to the
following remarkable integrable Dubrovin equation for the function
$f (u^2, u^3)$: \be {\pa^3 f \over \pa (u^3)^3} = \left ( {\pa^3 f
\over \pa (u^2)^2 \pa u^3} \right )^2 - {\pa^3 f \over \pa (u^2)^3}
{\pa^3 f \over \pa u^2 \pa (u^3)^2}. \label{f} \ee We introduce here
new independent variables: $x = u^2$, $t = u^3$. The equation
(\ref{f}) takes the form \be f_{ttt} =  (f_{xxt})^2 - f_{xxx}
f_{xtt}. \label{f1} \ee This equation is connected to quantum
cohomology of projective plane and classical problems of enumerative
geometry (see [17]).

It was proved by the present author in [18] (see also [19],
[10]--[12]) that the equation (\ref{f1}) is equivalent to the
integrable non-diagonalizable system of hydrodynamic type \be \left
(
\begin{array} {c} a^1\\ a^2\\ a^3
\end{array} \right )_t =
\left ( \begin{array} {ccc} 0 & 1 & 0\\  0 & 0 & 1\\
- a^3 & 2 a^2 & - a^1
\end{array} \right )  \left ( \begin{array} {c} a^1\\ a^2\\ a^3
\end{array} \right )_x, \label{syst}
\ee \be a^1 = f_{xxx}, \ \ \ a^2 = f_{xxt},\ \ \ a^3 = f_{xtt}. \ee

The first Hamiltonian structure of system (\ref{syst}) given by a
Dubrovin--Novikov bracket was found in [20]: \be \{ I, J\}_1 = \int
{\delta I \over \delta a^i (x)} M_1^{ij} {\delta J \over \delta a^j
(x)} dx, \label{ham1}\ee \bea  && M_1 = (M_1^{ij}) =
\left ( \begin{array} {ccc} - {3 \over 2} & {1 \over 2} a^1 & a^2\\
{1 \over 2} a^1 & a^2 & {3 \over 2} a^3\\ a^2 & {3 \over 2} a^3 & 2
((a^2)^2 - a^1 a^3)
\end{array} \right ) {d \over dx} + \nn\\ && +
\left ( \begin{array} {ccc} 0 & {1 \over 2} a^1_x & a^2_x\\
 0 & {1 \over 2} a^2_x & a^3_x\\
0 & {1 \over 2} a^3_x & ((a^2)^2 - a^1 a^3)_x
\end{array} \right ). \label{ham1a}
\eea The metric \be  (g_1^{ij} (a)) =
\left ( \begin{array} {ccc} - {3 \over 2} & {1 \over 2} a^1 & a^2\\
{1 \over 2} a^1 & a^2 & {3 \over 2} a^3\\ a^2 & {3 \over 2} a^3 & 2
((a^2)^2 - a^1 a^3)
\end{array} \right ) \label{met1}
\ee is flat and the Poisson bracket of hydrodynamic type
(\ref{ham1}), (\ref{ham1a}) is local (a Dubrovin--Novikov bracket).
The functional \be H_1 = \int a^3 dx \ee is the corresponding
Hamiltonian of system (\ref{syst}).

The bi-Hamiltonian structure of system (\ref{syst}) was found in
[10] (see also [11], [12]). The second Hamiltonian structure of
system (\ref{syst}) is given by a homogeneous third-order
Dubrovin--Novikov bracket: \be \{ I, J\}_2 = \int {\delta I \over
\delta a^i (x)} M_2^{ij} {\delta J \over \delta a^j (x)} dx,
\label{ham2}\ee \bea && M_2 = (M_2^{ij}) =
\left ( \begin{array} {ccc} 0 & 0 & 1\\
0 & 1 & - a^1\\ 1 & - a^1 & (a^1)^2 + 2 a^2
\end{array} \right ) \left ({d \over dx} \right )^3 + \nn\\ && +
\left ( \begin{array} {ccc} 0 & 0 & 0\\
0 & 0 & - 2 a^1_x\\ 0 & - a^1_x & 3 (a^2_x + a^1 a^1_x)
\end{array} \right ) \left ({d \over dx} \right )^2 + \nn\\ && +
\left ( \begin{array} {ccc} 0 & 0 & 0\\
 0 & 0 & 0\\
0 & 0 & a^2_{xx} + (a^1_x)^2 + a^1  a^1_{xx}
\end{array} \right ) {d \over dx}. \label{ham2a}
\eea The second Poisson bracket (\ref{ham2}), (\ref{ham2a}) is
compatible with the first Poisson bracket (\ref{ham1}),
(\ref{ham1a}).

The metric \be  (g_1^{ij} (a)) = \left ( \begin{array} {ccc} 0 & 0 & 1\\
0 & 1 & - a^1\\ 1 & - a^1 & (a^1)^2 + 2 a^2
\end{array} \right )
 \label{met2}
\ee is flat. The non-local functional \be H_2 = - \int \left ( {1
\over 2} a^1 \left ( \left ( {d \over dx} \right )^{- 1} a^2 \right
)^2 + \left ( \left ( {d \over dx} \right )^{- 1} a^2 \right ) \left
( \left ( {d \over dx} \right )^{- 1} a^3 \right ) \right ) dx \ee
is the corresponding Hamiltonian of system (\ref{syst}).

First of all, we note that the class of non-locally bi-Hamiltonian
systems of hydrodynamic type (\ref{nonlmS2})--(\ref{nonlm5}) is very
rich, there are many well-known important examples arising in
various applications. An explicit general construction of locally
and non-locally bi-Hamiltonian systems of hydrodynamic type and
corresponding integrable hierarchies that are generated by pairs of
compatible Poisson brackets of hydrodynamic type and integrable
description of local and non-local compatible Poisson brackets of
hydrodynamic type were found and studied by the present author in
[21]--[33], [4] (see also [34], [35], [13]).

\begin{center}
{\large \bf {Compatible metrics and non-locally bi-Hamiltonian
systems of hydrodynamic type}}
\end{center}

Now recall some necessary basic facts of the general theory of
compatible metrics [27]--[32]. Two Riemannian or pseudo-Riemannian
contravariant metrics $g_1^{ij} (u)$ and $g_2^{ij} (u)$ are called
{\it compatible} if for any linear combination of these metrics \be
g^{ij} (u) = \lambda_1 g_1^{ij} (u) + \lambda_2 g_2^{ij} (u),
\label{comb} \ee where $\lambda_1$ and $\lambda_2$ are arbitrary
constants such that $\det ( g^{ij} (u) ) \neq 0$, the coefficients
of the corresponding Levi--Civita connections and the components of
the corresponding Riemannian curvature tensors are related by the
same linear formula [27]--[29]: \be \Gamma^{ij}_k (u) = \lambda_1
\Gamma^{ij}_{1, k} (u) + \lambda_2 \Gamma^{ij}_{2, k} (u),
\label{sv} \ee \be R^{ij}_{kl} (u) = \lambda_1 R^{ij}_{1, kl} (u) +
\lambda_2 R^{ij}_{2, kl} (u). \label{kr} \ee The indices of the
coefficients of the Levi--Civita connections
 $\Gamma^i_{jk} (u)$
and the indices of the Riemannian curvature tensors $R^i_{jkl} (u)$
are raised and lowered by the metrics corresponding to them: \bea &&
\Gamma^{ij}_k (u) = g^{is} (u) \Gamma^j_{sk} (u),
 \ \ \ \Gamma^i_{jk} (u) = {1 \over 2} g^{is} (u) \left (
{\pa g_{sk} \over \pa u^j} + {\pa g_{js} \over \pa u^k} -
{\pa g_{jk} \over \pa u^s} \right ),\nn\\
&&
R^{ij}_{kl} (u) = g^{is} (u) R^j_{skl} (u), \nn\\
&& R^i_{jkl} (u) = {\pa \Gamma^i_{jl} \over \pa u^k} - {\pa
\Gamma^i_{jk} \over \pa u^l} + \Gamma^i_{pk} (u) \Gamma^p_{jl} (u) -
\Gamma^i_{pl} (u) \Gamma^p_{jk} (u).\nn \eea

Two Riemannian or pseudo-Riemannian contravariant metrics $g_1^{ij}
(u)$ and $g_2^{ij} (u)$ are called {\it almost compatible} if for
any linear combination of these metrics (\ref{comb}) relation
(\ref{sv}) holds [27]--[29].

Let us introduce the affinor \be v^i_j (u) = g_1^{is} (u) g_{2, sj}
(u) \label{aff} \ee and consider the Nijenhuis tensor of this
affinor \be N^k_{ij} (u) = v^s_i (u) {\pa v^k_j \over \pa u^s} -
v^s_j (u) {\pa v^k_i \over \pa u^s}  + v^k_s (u) {\pa v^s_i \over
\pa u^j} - v^k_s (u) {\pa v^s_j \over \pa u^i}.  \label{nij} \ee

{\bf Theorem 1 [27]--[29].} {\it Any two metrics $g_1^{ij} (u)$ and
$g_2^{ij} (u)$ are almost compatible if and only if the
corresponding Nijenhuis tensor $N^k_{ij} (u)$ {\rm (\ref{nij})}
vanishes.}

Assume that a pair of metrics $g^{ij}_1 (u)$ and $g^{ij}_2 (u)$ is
non-singular, i.e., the eigenvalues of this pair of metrics are
distinct. Furthermore, assume that the metrics $g_1^{ij} (u)$ and
$g_2^{ij} (u)$ are almost compatible, i.e., the corresponding
Nijenhuis tensor $N^k_{ij} (u)$ (\ref{nij}) vanishes. It was proved
in our papers [27]--[29] that, in this case, the metrics $g^{ij}_1
(u)$ and $g^{ij}_2 (u)$ are compatible, i.e., relation (\ref{kr})
holds.

It is obvious that the eigenvalues of the pair of metrics $g^{ij}_1
(u)$ and $g^{ij}_2 (u)$  coincide with the eigenvalues of the
affinor $v^i_j (u)$ (\ref{aff}). But it is well known that if all
eigenvalues of an affinor are distinct, then it always follows from
the vanishing of the Nijenhuis tensor of this affinor that there
exist special local coordinates ({\it Riemann invariants}) such
that, in these coordinates, the affinor reduces to a diagonal form
in the corresponding neighbourhood [7] (see also [6]).

Hence, we can consider that the affinor $v^i_j (u)$ is diagonal in
the local coordinates (Riemann invariants) $u^1,...,u^N$, i.e., \be
 v^i_j (u) = f^i (u) \delta^i_j,
\ee where is no summation over the index $i$. By our assumption, the
eigenvalues $f^i (u), \ i=1,...,N,$ coinciding with the eigenvalues
of the pair of metrics $g^{ij}_1 (u)$ and $ g^{ij}_2 (u)$ are
distinct: \be f^i (u) \neq f^j (u) \ \ {\rm \ if\ } \ \ i\neq j. \ee

{\bf Lemma 1.}  {\it If the affinor $v^i_j (u)$ {\rm (\ref{aff})} is
diagonal in certain local coordinates {\rm(}Riemann invariants{\rm)}
and all its eigenvalues are distinct, then, in these coordinates,
the metrics $g^{ij}_1 (u) $ and $g^{ij}_2 (u)$ are also necessarily
diagonal, i.e., in this case both the metrics $g^{ij}_1 (u) $ and
$g^{ij}_2 (u)$ are diagonal in the Riemann invariants.}

Actually, we have
$$g^{ij}_1 (u) = f^i (u) g^{ij}_2 (u).$$
It follows from the symmetry of the metrics $g^{ij}_1 (u)$ and $
g^{ij}_2 (u)$ that for any indices $i$ and  $j$ \be (f^i (u) - f^j
(u))  g^{ij}_2 (u) = 0, \ee where is no summation over indices,
i.e.,
$$ g^{ij}_2 (u) = g^{ij}_1 (u) =0 \ \
{\rm \  if \ } \ \ i\neq  j. $$

{\bf Lemma 2.}  {\it Let an affinor $w^i_j (u)$ be diagonal in
certain local coordinates {\rm(}Riemann invariants{\rm)} $u=
(u^1,...,u^N)$, i.e., $w^i_j (u) =\mu^i (u) \delta^i_j$. \bitem
\item[1)] If all the eigenvalues $\mu^i (u), \ i=1,...,N,$
of the diagonal affinor are distinct, i.e., $\mu^i (u) \neq \mu^j
(u)$ for $i \neq  j$, then the Nijenhuis tensor of this affinor
vanishes if and only if the $i$th eigenvalue $\mu^i (u)$ depends
only on the coordinate $u^i.$
\item[2)] If all the eigenvalues coincide, then
the Nijenhuis tensor vanishes.
\item[3)] In the general case of an arbitrary diagonal affinor
$w^i_j (u) =\mu^i (u) \delta^i_j$, the Nijenhuis tensor vanishes if
and only if \be {\pa \mu^i \over \pa u^j} = 0 \ee for all indices
$i$ and $j$ such that $\mu^i (u) \neq \mu^j (u).$ \eitem }

It follows from Lemmas 1 and 2 that for any non-singular pair of
almost compatible metrics there always exist local coordinates
(Riemann invariants) in which the metrics have the form
$$g^{ij}_2 (u) = g^i (u) \delta^{ij}, \ \ \
g^{ij}_1 (u) = f^i (u^i) g^i (u) \delta^{ij}.$$

Moreover, any pair of diagonal metrics of the form $g^{ij}_2 (u) =
g^i (u) \delta^{ij}$ and $g^{ij}_1 (u) = f^i (u^i) g^i (u)
\delta^{ij}$ for any nonzero functions $f^i (u^i),$ $ i=1,...,N,$
(here they can be, for example, coinciding nonzero constants, i.e.,
the pair of metrics may be ``singular'') is almost compatible, since
the corresponding Nijenhuis tensor always vanishes for any pair of
metrics of this form. It was proved in our papers [28], [29] that an
arbitrary pair of diagonal metrics of such the form, $g^{ij}_2 (u) =
g^i (u) \delta^{ij}$ and $g^{ij}_1 (u) = f^i (u^i) g^i (u)
\delta^{ij}$, for arbitrary nonzero functions $f^i (u^i),$ $
i=1,...,N,$ (the pair of metrics may be ``singular''), is always
compatible, i.e., in this case, relation (\ref{kr}) holds. We note
that, as it was shown in [27]--[29], in general almost compatible
metrics are not necessarily compatible even in the case of flat
metrics or metrics of constant curvature, i.e., in the case of the
Dubrovin--Novikov or the Mokhov--Ferapontov brackets, but if a pair
of almost compatible metrics is not compatible, then this pair of
metrics must be singular. Thus, we have the following important
statements.

{\bf Theorem 2 [27]--[29].} {\it If a pair of metrics $g_1^{ij} (u)$
and $g_2^{ij} (u)$ is non-singular, i.e., the roots of the equation
\be \det ( g_1^{ij} (u) -  \lambda g_2^{ij} (u)) =0 \ee are
distinct, then it follows from the vanishing of the Nijenhuis tensor
of the affinor $v^i_j (u) = g_1^{is} (u) g_{2, sj} (u)$ that the
metrics $g_1^{ij} (u)$ and $g_2^{ij} (u)$ are compatible. Thus, a
non-singular pair of metrics is compatible if and only if the
metrics are almost compatible.}

{\bf Theorem 3 [27]--[29].} {\it An arbitrary non-singular pair of
metrics is compatible if and only if there exist local coordinates
{\rm(}Riemann invariants{\rm)} $u = (u^1,...,u^N)$ such that both
the metrics are diagonal in these coordinates and have the following
special form: $g^{ij}_2 (u) = g^i (u) \delta^{ij}$ and $g^{ij}_1 (u)
= f^i (u^i) g^i (u) \delta^{ij},$ where one of the metrics, here
$g^{ij}_2 (u),$ is an arbitrary diagonal metric and $f^i (u^i),$
$i=1,...,N,$ are arbitrary {\rm(}generally speaking, complex{\rm)}
nonzero functions of single variables. If some of the functions $f^i
(u^i),$ $i=1,...,N,$ are coinciding nonzero constants, then the pair
of metrics of this form is singular but, nevertheless, compatible.}

{\bf Theorem 4 [28].} {\it If non-local Poisson brackets of
hydrodynamic type {\rm (\ref{nonlm3a})--(\ref{nonlm5})} are
compatible, then their metrics are compatible.}

In [2] Ferapontov proved that a bracket
(\ref{nonlm1})--(\ref{nonlm3}) is a Poisson bracket, i.e., it is
skew-symmetric and satisfies the Jacobi identity, if and only if
\bitem
\item [(1)] $b^{ij}_k (u) = - g^{is} (u) \Gamma ^j_{sk} (u),$ where
$\Gamma^j_{sk} (u)$ is the Riemannian connection generated by the
contravariant metric $g^{ij} (u)$
 (the Levi--Civita connection),

\item [(2)] the pseudo-Riemannian
metric $g^{ij} (u)$ and the set of affinors $(w_n)^i_j (u)$ satisfy
the relations:

\be g_{ik} (u) (w_n)^k_j (u) = g_{jk} (u) (w_n)^k_i (u), \ \ \ n =
1,...,L,  \label{peter1} \ee \be \nabla_k (w_n)^i_j (u) = \nabla_j
(w_n)^i_k (u), \ \ \ n = 1,...,L, \label{peter2} \ee \be R^{ij}_{kl}
(u) =  \sum_{m =1}^L \sum_{n =1}^L \mu^{mn} \left ( (w_m)^i_l (u)
(w_n)^j_k (u) - (w_m)^j_l (u) (w_n)^i_k (u) \right ). \label{gauss}
\ee Moreover, the family of affinors $w_n (u)$ is commutative:
$[w_m, w_n] =0.$ \eitem

 If non-local Poisson brackets of hydrodynamic type
(\ref{nonlm3a})--(\ref{nonlm5}) are compatible, then it follows from
the conditions of compatibility and from Ferapontov's theorem that,
first, relation (\ref{sv}) holds, i.e., the metrics $g_1^{ij} (u)$
and $g_2^{ij} (u)$ are almost compatible, and, secondly, the
curvature tensor for the metric $g^{ij} (u) = \lambda_1 g^{ij}_1 (u)
+ \lambda_2 g^{ij}_2 (u)$ has the form \bea && R^{ij}_{kl} (u) =
\sum_{m =1}^{L_1} \sum_{n =1}^{L_1} \lambda_1 \mu_1^{mn} \left (
(w_{1, m})^i_l (u) (w_{1, n})^j_k (u)
- (w_{1, m})^j_l (u) (w_{1, n})^i_k (u) \right ) + \nn\\
&& + \sum_{m =1}^{L_2} \sum_{n =1}^{L_2} \lambda_2 \mu_2^{mn} \left
( (w_{2, m})^i_l (u) (w_{2, n})^j_k (u) - (w_{2, m})^j_l (u) (w_{2,
n})^i_k (u) \right ) = \nn\\ && = \lambda_1 R^{ij}_{1, kl} (u) +
\lambda_2 R^{ij}_{2, kl} (u),\nn \eea i.e., relation (\ref{kr})
holds and hence the metrics $g_1^{ij} (u)$ and $g_2^{ij} (u)$ are
compatible.

{\bf Theorem 5 [28].} {\it Let two non-local Poisson brackets of
hydrodynamic type {\rm(\ref{nonlm3a})--(\ref{nonlm5})} correspond to
submanifolds with holonomic net of curvature lines and be given in
coordinates of curvature lines. In this case, if the corresponding
pair of metrics is non-singular, then the non-local Poisson brackets
of hydrodynamic type are compatible if and only if their metrics are
compatible.}

In this case the metrics $g^{ij}_1 (u) = g^i_1 (u) \delta^{ij}$ and
$g^{ij}_2 (u) = g^i_2 (u) \delta^{ij}$, and also the Weingarten
operators $(w_{1, n})^i_j (u) = (w_{1, n})^i (u) \delta^i_j$ and
$(w_{2, n})^i_j (u) = (w_{2, n})^i (u) \delta^i_j$  are diagonal in
the coordinates under consideration. For any such ``diagonal'' case,
condition (\ref{peter1}) is automatically fulfilled, all the
Weingarten operators commute, conditions (\ref{peter2}) and
(\ref{gauss}) have the following form, respectively: \be 2 g^i (u)
{\pa (w_n)^i \over \pa u^k} = ((w_n)^i - (w_n)^k) {\pa g^i \over \pa
u^k} \ \ \ \ {\rm \ for \ all \ } \ \ i \neq k,\label{hol1} \ee \bea
&& R^{ij}_{j\, i} (u) =  \sum_{m =1}^L \sum_{n =1}^L \mu^{mn}
(w_m)^i (u) (w_n)^j (u), \nn\\ && R^{ij}_{kl} (u) = 0  \ \ {\rm \ if
\ } \ i \neq k, \ i \neq l, \ \ {\rm \ or \ if \ } \ j \neq k, \ j
\neq l. \label{hol2} \eea

It follows from non-singularity of the pair of the metrics and from
compatibility of the metrics that the corresponding Nijenhuis tensor
vanishes and there exist functions
 $f^i (u^i), i=1,...,N,$ such that:
$$g^i_1 (u) = f^i (u^i) g^i_2 (u).$$
Using relations (\ref{hol1}) and (\ref{hol2}), it is easy to prove
that in this case it follows from compatibility of the metrics that
an arbitrary linear combination of non-local Poisson brackets under
consideration is also a Poisson bracket.

{\bf Theorem 5 [33].} {\it If the pair of metrics $g^{ij}_1 (u)$ and
$g^{ij}_2 (u)$ is non-singular, then the non-local Poisson brackets
of hydrodynamic type $\{I, J \}_1$ and $\{ I, J \}_2$
{\rm(\ref{nonlm3a})--(\ref{nonlm5})} are compatible if and only if
the metrics are compatible and both the metrics $g^{ij}_1 (u),$
$g^{ij}_2 (u)$ and the affinors $(w_{1, n})^i_j (u),$ $(w_{2,
n})^i_j (u)$ can be samultaneously diagonalized in a domain of local
coordinates.}

It is sufficient to prove here that if the pair of metrics is
non-singular and the Poisson brackets are compatible, then both the
metrics $g^{ij}_1 (u),$ $g^{ij}_2 (u)$ and the affinors $(w_{1,
n})^i_j (u),$ $(w_{2, n})^i_j (u)$ can be samultaneously
diagonalized in a domain of local coordinates. All the rest was
already proved above. First of all, it was proved that in this case
the metrics $g^{ij}_1 (u)$ and $g^{ij}_2 (u)$ are compatible. Since
the pair of metrics is non-singular, there exist local coordinates
such that the metrics are diagonal and have the following special
form in these coordinates: $g^{ij}_2 (u) = g^i (u) \delta^{ij}$ and
$g^{ij}_1 (u) = f^i (u^i) g^i (u) \delta^{ij},$ where $f^i (u^i),$
$1 \leq i \leq N,$ are functions of single variable. The functions
$f^i (u^i)$ are the eigenvalues of the pair of metrics $g^{ij}_1
(u)$ and $g^{ij}_2 (u)$, therefore they are distinct by assumption
of the theorem even in the case if they are constants (they can not
be coinciding constants). It follows from the compatibility of the
Poisson brackets $\{ I, J \}_1$ and $\{ I, J \}_2$ (it is necessary
to consider relation (\ref{04}) for the pencil $\{ I, J \}_1 +
\lambda \{ I, J \}_2$) that \be g^{is}_1 (w_{2, n})^j_s = g^{js}_1
(w_{2, n})^i_s, \label{k1} \ee \be g^{is}_2 (w_{1, n})^j_s =
g^{js}_2 (w_{1, n})^i_s. \label{k2} \ee Besides, from relation
(\ref{04}) for the Poisson brackets $\{ I, J \}_1$ and $\{ I, J
\}_2$ we have \be g^{is}_1 (w_{1, n})^j_s = g^{js}_1 (w_{1, n})^i_s,
\label{k3} \ee \be g^{is}_2 (w_{2, n})^j_s = g^{js}_2 (w_{2,
n})^i_s. \label{k4} \ee From (\ref{k1}) and (\ref{k4}) in our
special local coordinates we obtain \be g^i (w_{2, n})^j_i = g^j
(w_{2, n})^i_j, \ee \be f^i (u^i) g^i (w_{2, n})^j_i = f^j (u^j) g^j
(w_{2, n})^i_j. \ee Therefore \be (w_{2, n})^i_j = {g^i \over g^j}
(w_{2, n})^j_i = {f^i (u^i) g^i \over f^j (u^j) g^j} (w_{2, n})^j_i,
\ee i.e., \be \left ( 1 - {f^i (u^i) \over f^j (u^j)} \right )
(w_{2, n})^j_i = 0. \ee Consequently, since all the functions $f^i
(u^i)$ are distinct, we get \be (w_{2, n})^j_i = 0 \ \ {\rm \ for \
} i \neq j. \ee Similarly, from (\ref{k2}) and (\ref{k3}) we have
\be (w_{1, n})^j_i = 0 \ \ {\rm \ for \  } i \neq j. \ee Thus, both
the metrics $g^{ij}_1 (u),$ $g^{ij}_2 (u)$ and the affinors $(w_{1,
n})^i_j (u),$ $(w_{2, n})^i_j (u)$ are diagonal in our special local
coordinates.

Now we can prove the main theorem of the paper.

{\bf Theorem 6.} {\it For an arbitrary non-singular
{\rm(}semisimple{\rm)} non-locally bi-Hamiltonian system of
hydrodynamic type {\rm(\ref{nonlmS2})--(\ref{nonlm5})}, there exist
local coordinates {\rm(}Riemann invariants{\rm)} such that all the
related matrix differential-geometric objects, namely, the matrix
$V^i_j (u)$ of this system of hydrodynamic type, the metrics
$g^{ij}_1 (u)$ and $g^{ij}_2 (u)$ and the affinors $(w_{1, n})^i_j
(u)$ and $(w_{2, n})^i_j (u)$ of the non-local bi-Hamiltonian
structure of this system, are diagonal in these local coordinates.}

If we have a non-singular (semisimple) non-locally bi-Hamiltonian
system of hydrodynamic type (\ref{nonlmS2})--(\ref{nonlm5}), then it
was proved above that the metrics $g^{ij}_1 (u)$ and $g^{ij}_2 (u)$
of the non-local bi-Hamiltonian structure of this system are
compatible and there exist local coordinates such that $g^{ij}_2 (u)
= g^i (u) \delta^{ij}$ and $g^{ij}_1 (u) = f^i (u^i) g^i (u)
\delta^{ij},$ where $f^i (u^i),$ $i=1,...,N,$ are distinct nonzero
functions of single variable (generally speaking, complex), $f^i
(u^i) \neq f^j (u^j)$, $i \neq j$. It was also proved above that the
affinors $(w_{1, n})^i_j (u)$ and $(w_{2, n})^i_j (u)$ of the
non-local bi-Hamiltonian structure of this system, are diagonal in
these local coordinates. Let us prove that the matrix $V^i_j (u)$ of
this system is also diagonal in these special local coordinates.

Indeed, in these local coordinates, we have from relations
(\ref{jj1}): \be g_i (u) V^i_j (u) = g_j (u) V^j_i (u),\ \ \ \ f_i
(u^i) g_i (u) V^i_j (u) = f_j (u^j) g_j (u) V^j_i (u).\ee

Hence, \be V^i_j (u) = {g_j (u) \over g_i (u)} V^j_i (u) = {f_j
(u^j) g_j (u) \over f_i (u^i) g_i (u)} V^j_i (u), \ee i.e., \be {g_j
(u) \over g_i (u)} V^j_i (u) = {f_j (u^j) g_j (u) \over f_i (u^i)
g_i (u)} V^j_i (u).\ee Thus, \be (f_i (u^i) - f_j (u^j)) V^j_i (u) =
0,\ee i.e., \be V^i_j (u) = 0, \ \ \ i \neq j,\ee and the
diagonalizability of an arbitrary non-singular (semisimple)
non-locally bi-Ha\-miltonian system of hydrodynamic type is proved:
\be V^i_j (u) = V^i (u) \delta^i_j.\ee

The diagonalizability of non-singular (semisimple) locally
bi-Hamiltonian systems of hydrodynamic type
(\ref{nonlmS2})--(\ref{nonlm5}) was noticed in [35]; it follows
immediately from the theory of non-singular pairs of compatible flat
metrics [27]--[34].

We note that it does not follow from the proof that $V^i (u) \neq
V^j (u)$ if $i \neq j$, i.e., an arbitrary non-singular (semisimple)
non-locally bi-Hamiltonian system of hydrodynamic type must not be
necessarily strictly hyperbolic but we do not know examples of such
systems with some coinciding eigenvalues (velocities) $V^i (u)$. We
conjecture that there exist such systems and this is a very
interesting problem to find non-singular (semisimple) non-locally
bi-Hamiltonian systems of hydrodynamic type that are not strictly
hyperbolic (i.e., they have some coinciding eigenvalues $V^i (u)$).

We also note that the non-singularity  condition of the pair of
metrics is very essential and we conjecture that there exist
non-diagonalizable singular non-locally bi-Hamiltonian systems of
hydrodynamic type. It is also an interesting problem to find
examples of non-diagonalizable singular non-locally bi-Hamiltonian
systems of hydrodynamic type.

\medskip

{\bf {Acknowledgements.}} The work was supported by the
Max-Planck-Institut f\"ur Mathematik (Bonn, Germany), by the Russian
Foundation for Basic Research (project no.~08-01-00464) and by a
grant of the President of the Russian Federation (project no.
NSh-1824.2008.1).

\begin{center}
\bf {References}
\end{center}

\medskip

[1] B. A. Dubrovin and S. P. Novikov, ``The Hamiltonian formalism of
one-dimensi\-onal systems of hydrodynamic type and the
Bogolyubov--Whitham averaging method,'' Dokl. Akad. Nauk SSSR, Vol.
270, No. 4, 1983, pp. 781--785; English translation in Soviet Math.
Dokl., Vol. 27, 1983, pp. 665--669.

[2] E. V. Ferapontov, ``Differential geometry of nonlocal
Hamiltonian operators of hydrodynamic type,'' Funkts. Analiz i Ego
Prilozh., Vol. 25, No. 3, 1991, 37--49; English translation in
Functional Analysis and its Applications, Vol. 25, No. 3, 1991, pp.
195--204.

[3] O. I. Mokhov and E. V. Ferapontov, ``Non-local Hamiltonian
operators of hydrodynamic type related to metrics of constant
curvature,'' Uspekhi Matem. Nauk, Vol. 45, No. 3, 1990, pp.
191--192; English translation in Russian Mathematical Surveys, Vol.
45, No. 3, 1990, pp. 218--219.

[4] O. I. Mokhov, ``Nonlocal Hamiltonian operators of hydrodynamic
type with flat metrics, integrable hierarchies, and the
associativity equations,'' Funkts. Analiz i Ego Prilozh., Vol. 40,
No. 1, 2006, pp. 14--29; English translation in Functional Analysis
and its Applications, Vol. 40, No. 1, 2006, pp. 11--23; \noindent
http://arXiv.org/math.DG/0406292 (2004).

[5] S. P. Tsarev, ``Geometry of Hamiltonian systems of hydrodynamic
type. The generalized hodograph method,'' Izvestiya Akad. Nauk SSSR,
Ser. Matem., Vol. 54, No. 5, 1990, pp. 1048--1068; English
translation in Math. USSR -- Izvestiya, Vol. 54, No. 5, 1990, pp.
397--419.

[6] A. Haantjes, ``On $X_{n-1}$-forming sets of eigenvectors,''
Indagationes Mathematicae, Vol. 17, No. 2, 1955, pp. 158--162.

[7] A. Nijenhuis, ``$X_{n-1}$-forming sets of eigenvectors,''
Indagationes Mathematicae, Vol. 13, No. 2, 1951, pp. 200--212.

[8] M. V. Pavlov, R. A. Sharipov and S. I. Svinolupov, ``Invariant
integrability criterion for equations of hydrodynamic type,''
Funkts. Analiz i Ego Prilozh., Vol. 30, No. 1, 1996, pp. 18--29;
English translation in Functional Analysis and its Applications,
Vol. 30, No. 1, 1996, pp. 15--22.

[9] F. Magri, ``A simple model of the integrable Hamiltonian
equation,'' J. Math. Phys., Vol. 19, No. 5, 1978, pp. 1156--1162.

[10] E. V. Ferapontov, C. A. P. Galv\~ao, O. I. Mokhov and Y. Nutku,
``Bi-Hamiltonian structure of equations of associativity in 2D
topological field theory,'' Comm. Math. Phys., Vol. 186, 1997, pp.
649--669.

[11] O. I. Mokhov, ``Symplectic and Poisson structures on loop
spaces of smooth manifolds, and integrable systems'', Uspekhi
Matematicheskikh Nauk, Vol. 53, No. 3, 1998, pp. 85--192; English
translation in Russian Mathematical Surveys, Vol. 53, No. 3, 1998,
pp. 515--622.

[12] O. I. Mokhov, ``Symplectic and Poisson geometry on loop spaces
of smooth manifolds and integrable equations'', Moscow--Izhevsk,
Institute of Computer Studies, 2004 (In Russian); English version:
Reviews in Mathematics and Mathematical Physics, Vol. 11, Part 2,
Harwood Academic Publishers, 2001; Second Edition: Reviews in
Mathematics and Mathematical Physics, Vol. 13, Part 2, Cambridge
Scientific Publishers, 2009.

[13] B. Dubrovin, ``Geometry of 2D topological field theories,'' In:
Integrable Systems and Quantum Groups, Lecture Notes in Math., Vol.
1620, Springer-Verlag, Berlin, 1996, pp. 120--348;
http://arXiv.org/hep-th/9407018 (1994).

[14] E. Witten, ``On the structure of the topological phase of
two-dimensional gravity,'' Nuclear Physics B, Vol. 340, 1990, pp.
281--332.

[15] E. Witten, ``Two-dimensional gravity and intersection theory on
moduli space,'' Surveys in Diff. Geometry, Vol. 1, 1991, pp.
243--310.

[16] R. Dijkgraaf, H. Verlinde and E. Verlinde, ``Topological
strings in $d < 1$,'' Nuclear Physics B, Vol. 352, 1991, pp. 59--86.

[17] M. Kontsevich and Yu. Manin, ``Gromov--Witten classes, quantum
cohomology, and enumerative geometry,'' Comm. Math. Phys., Vol. 164,
1994, pp. 525--562.

[18] O. I. Mokhov, ``Symplectic and Poisson geometry on loop spaces
of manifolds and nonlinear equations,'' In: Topics in Topology and
Mathematical Physics, Ed. S.P.Novikov, Amer. Math. Soc., Providence,
RI, 1995, pp. 121--151;

\noindent  http://arXiv.org/hep-th/9503076 (1995).

[19] O. I. Mokhov, ``Poisson and symplectic geometry on loop spaces
of smooth manifolds,'' In: Geometry from the Pacific Rim,
Proceedings of the Pacific Rim Geometry Conference held at National
University of Singapore, Republic of Singapore, December 12--17,
Eds. A.J.Berrick, B.Loo, H.-Y.Wang, Walter de Gruyter, Berlin, 1997,
pp. 285--309.

[20] O. I. Mokhov and E. V. Ferapontov, ``Equations of associativity
in two-dimen\-sional topological field theory as integrable
Hamiltonian non-diagonalizable systems of hydrodynamic type,''
Funkts. Analiz i Ego Prilozh., Vol. 30, No. 3, 1996, pp. 62--72;
English translation in Functional Analysis and its Applications,
Vol. 30, No. 3, 1996, pp. 195--203; \noindent
http://arXiv.org/hep-th/9505180 (1995).

[21] O. I. Mokhov, ``Compatible nonlocal Poisson brackets of
hydrodynamic type and integrable hierarchies related to them,''
Teoret. Matem. Fizika, Vol. 132, No. 1, 2002, pp. 60--73; English
translation in Theoretical and Mathematical Physics, Vol. 132, No.
1, 2002, pp. 942--954; \noindent http://arXiv.org/math.DG/0201242
(2002).

[22] O. I. Mokhov, ``Integrable bi-Hamiltonian systems of
hydrodynamic type,'' Uspekhi Matematicheskikh Nauk, Vol. 57, No. 1,
2002, pp. 157--158; English translation in Russian Mathematical
Surveys, Vol. 57, No. 1, 2002, pp. 153--154.

[23] O. I. Mokhov, ``Integrable bi-Hamiltonian hierarchies generated
by compatible metrics of constant Riemannian curvature,'' Uspekhi
Matematicheskikh Nauk, Vol. 57, No. 5, 2002, pp. 157--158; English
translation in Russian Mathematical Surveys, Vol. 57, No. 5, 2002,
pp. 999--1001.

[24] O. I. Mokhov, ``The Liouville canonical form for compatible
nonlocal Poisson brackets of hydrodynamic type and integrable
hierarchies,'' Funkts. Analiz i Ego Prilozh., Vol. 37, No. 2, 2003,
pp. 28--40; English translation in Functional Analysis and its
Applications, Vol. 37, No. 2, 2003, pp. 103--113;

\noindent http://arXiv.org/math.DG/0201223(2002).

[25] O. I. Mokhov, ``Frobenius manifolds as a special class of
submanifolds in pseudo-Euclidean spaces,''  In: Geometry, Topology,
and Mathematical Physics (S.P.Novikov's Seminar: 2006-2007), Amer.
Math. Soc., Providence, RI, 2008, pp. 213--246; arXiv: 0710.5860
(2007).

[26] O. I. Mokhov, ``Compatible Dubrovin--Novikov Hamiltonian
operators, Lie derivative, and integrable systems of Hydrodynamic
type,'' Teoret. Matem. Fizika, Vol. 133, No. 2, 2002, pp. 279--288;
English translation in Theoretical and Mathematical Physics, Vol.
133, No. 2, 2002, pp. 1557--1564; http://arXiv.org/math.DG/0201281
(2002).

[27] O. I. Mokhov, ``Compatible and almost compatible metrics,''
Uspekhi Matematicheskikh Nauk, Vol. 55, No. 4, 2000, pp. 217--218;
English translation in Russian Mathematical Surveys, Vol. 55, No. 4,
2000, pp. 819--821.

[28] O. I. Mokhov, ``Compatible and almost compatible
pseudo-Riemannian metrics,'' Funkts. Analiz i Ego Prilozh., Vol. 35,
No. 2, 2001, pp. 24--36; English translation in Functional Analysis
and its Applications, Vol. 35, No. 2, 2001, pp. 100--110;
http://arXiv.org/math.DG/0005051 (2000).

[29] O. I. Mokhov, ``Compatible flat metrics,'' Journal of Applied
Mathematics, Vol. 2, No. 7, 2002, pp. 337--370;
http://arXiv.org/math.DG/0201224 (2002).

[30] O. I. Mokhov, ``Flat pencils of metrics and integrable
reductions of the Lam\'e equations,'' Uspekhi Matematicheskikh Nauk,
Vol. 56, No. 2, 2001, pp. 221--222; English translation in Russian
Mathematical Surveys, Vol. 56, No. 2, 2001, pp. 416--418.

[31] O. I. Mokhov, ``Integrability of the equations for nonsingular
pairs of compatible flat metrics,'' Teoret. Matem. Fizika, Vol. 130,
No. 2, 2002, pp. 233--250; English translation in Theoretical and
Mathematical Physics, Vol. 130, No. 2, 2002, pp. 198--212;
http://arXiv.org/math.DG/0005081 (2000).

[32] O. I. Mokhov, ``Compatible metrics of constant Riemannian
curvature: local geometry, nonlinear equations, and integrability,''
Funkts. Analiz i Ego Prilozh., Vol. 36, No. 3, 2002, pp. 36--47;
English translation in Functional Analysis and its Applications,
Vol. 36, No. 3, 2002, pp. 196--204; http://arXiv.org/math.DG/0201280
(2002).

[33] O. I. Mokhov, ``Lax pairs for equations describing compatible
nonlocal Poisson brackets of hydrodynamic type and integrable
reductions of the Lame equations,'' Teoret. Matem. Fizika, Vol. 138,
No. 2, 2004, pp. 283--296; English translation in Theoretical and
Mathematical Physics, Vol. 138, No. 2, 2004, pp. 238--249;

\noindent http://arXiv.org/math.DG/0202036 (2002).

[34] E. V. Ferapontov, ``Compatible Poisson brackets of hydrodynamic
type,'' J. Phys. A: Math. Gen., Vol. 34, 2001, pp. 2377--3388;
http://arXiv.org/math.DG/0005221 (2000).

[35] B. Dubrovin, S.-Q. Liu and Y. Zhang, ``On Hamiltonian
perturbations of hyperbolic systems of conservation laws, I:
Quasi-triviality of bi- Hamiltonian perturbations,'' Comm. Pure
Appl. Math., Vol. 59, 2006, 559 - 615; arXiv:math/0410027 (2004).

\begin{flushleft}
{\bf O. I. Mokhov}\\
Centre for Nonlinear Studies,\\
L.D.Landau Institute for Theoretical Physics,\\
Russian Academy of Sciences,\\
Kosygina str., 2,\\
Moscow, 117940, Russia;\\
Department of Geometry and Topology,\\
Faculty of Mechanics and Mathematics,\\
M.V.Lomonosov Moscow State University,\\
Moscow, 119992, Russia\\
{\it E-mail\,}: mokhov@mi.ras.ru; mokhov@landau.ac.ru; mokhov@bk.ru\\
\end{flushleft}

\end{document}